\documentclass[11pt]{article}
\usepackage{amsmath,amsthm,amsfonts,latexsym,amssymb,enumerate,color}
\usepackage{graphicx}

\def\gA{\mathcal{A}}

\def\CC{\mathbb{C}}
\def\DD{\mathbb{D}}

\def\RR{\mathbb{R}}
\def\TT{\mathbb{T}}
\def\ZZ{\mathbb{Z}}

\def \codim{\mathop{\rm codim}\nolimits}
\def\dim{\mathop{\rm dim}\nolimits}
\def\GCD{\mathop{\rm GCD}\nolimits}
\def\ker{\mathop{\rm ker}\nolimits}

\def\Lat{\mathop{\rm Lat}\nolimits}

\renewcommand\span{\mathop{\rm span}\nolimits}
\def\spam{\span} % I always mistype this

\def\ol{\overline}

\renewcommand\phi{\varphi}

\newcommand{\simast}{%
    \ensuremath{%
       \stackrel{\mathsf{\ast}}{\sim}}}

\newtheorem{thm}{Theorem}[section]

\newtheorem{defn}[thm]{Definition}
\newtheorem{lem}[thm]{Lemma}
\newtheorem{cor}[thm]{Corollary}

\numberwithin{equation}{section}

\def\beginpf{\begin{proof}}
\def\endpf{\end{proof}}
\def\beq{\begin{equation}}
\def\eeq{\end{equation}}

\def\imag{\mathop{\rm Im}\nolimits}

\begin{document}

\title{Asymmetric truncated Toeplitz operators and Toeplitz operators with matrix symbol}

\date {}

\author{M.~Cristina C\^amara\thanks{
Center for Mathematical Analysis, Geometry, and Dynamical Systems,
Departamento de Matem\'{a}tica, Instituto Superior T\'ecnico, 1049-001 Lisboa, Portugal. \tt ccamara@math.ist.utl.pt} \ and  Jonathan R.~Partington\thanks{School of Mathematics, University of Leeds, Leeds LS2~9JT, U.K. \tt j.r.partington@leeds.ac.uk}
}

\maketitle

\begin{abstract}
Truncated Toeplitz operators and their asymmetric versions are studied in the context of the Hardy space
$H^p$ of the half-plane for $1<p<\infty$. It is shown that they are equivalent after extension to
$2 \times 2$ matricial Toeplitz operators, which allows one to deduce information about their
invertibility
properties. Shifted model spaces are presented in the context of invariant subspaces, allowing
one to deduce new Beurling--Lax theorems.
\end{abstract}

\noindent {\bf Keywords:}
Truncated Toeplitz operator, Toeplitz operator, equivalence by extension, model space, invariant subspace. 

\noindent{\bf MSC:} 47B35, 30H10.

%%%%%%%%%%%%%%%%%%%%%%%%%%%%%%%%%%%%%%%%%%%%%%%%%%%%%%%%%%%%

\section{Introduction}
\label{sec:1}

Truncated Toeplitz operators (TTO), a natural generalisation of finite Toe\-plitz matrices, have received
much attention since they were introduced by Sarason~\cite{sarason07}: see, for example
\cite{BCFT} and the recent survey \cite{GR13}. They appear in various contexts, for example in the
study of
finite Toeplitz matrices and finite-time convolution operators.
Here we treat a slightly more general class of operators, known as asymmetric
truncated Toeplitz operators (ATTO), although most of the results we prove are new even for
``standard'' TTO in the Hardy spaces $H^p$.\\

Full definitions and notation will be given later, but we work mostly with the Hardy spaces
$H_p^\pm$ of the upper and lower half-planes, for $1<p<\infty$, recalling the decomposition
$L_p(\RR)=H_p^+ \oplus H_p^-$. Many of our results may be rewritten for the disc, as we
shall see later, although the results often appear more complicated in this context.
 For an inner function $\theta\in H_\infty^+$
the {\em model space\/} $K_\theta$ may be defined as 
\beq\label{eq:defktheta}
K_\theta=H_p^+ \cap \theta H_p^-.
\eeq
 We then have
\beq\label{eq:decomp}
L_p(\RR)=H_p^- \oplus K_\theta \oplus \theta H_p^+,
\eeq
and we write $P_\theta$ to denote the associated projection $P_\theta: L_p(\RR) \to K_\theta$.

Then for $g \in L_\infty(\RR)$ the standard TTO $A_g^\theta$ is defined as follows.

\beq \label{II.1}
 A_g^\theta : K_\theta\rightarrow K_\theta, \quad  A_g^\theta=P_\theta(gI)_{\mid_{K_\theta}}=P_\theta(gI)_{\mid_{P_\theta L_p}}.
\eeq

If $\alpha$ and $\theta$ are inner functions, % that divides $\theta$ in $H_\infty^+$ (we write this $\alpha \preceq \theta$),
%then $K_\alpha \subseteq K_\theta$ and 
we define the operator $ A_g^{\alpha, \theta}$ as 
\beq 
 A_g^{\alpha, \theta} := P_\alpha(gI)_{\mid_{K_\theta}}=P_\alpha(gI)_{\mid_{P_\theta L_p}}.
\eeq

If $\alpha$ is an inner function that divides $\theta$ in $H_\infty^+$ (we write this $\alpha \preceq \theta$), let $P_{\alpha, \theta}$ denote $P_\theta-P_\alpha$, a projection with range
equal to the {\em shifted model space}
$K_{\alpha,\theta}:=\alpha K_{\overline\alpha \theta}$. Then we can define

 \beq 
 B_g^{\alpha, \theta} := P_{\alpha, \theta}(gI)_{\mid_{K_\theta}}=P_{\alpha, \theta}(gI)_{\mid_{P_\theta L_p}}.
\eeq

\vspace{0,1cm}

The operators $A_g^{\alpha, \theta}$ and $ B_g^{\alpha, \theta}$
 are particular cases of \emph{general WH operators} (see \cite{Speck}) in $L_p$, of the form
 \beq \label{II.4}
 P_1 A_{\mid_{P_2 L_p}}\,,
\eeq
where $P_1$ and $P_2$ are projections and $A$ is an operator in $L_p$. We say that $A_g^{\alpha, \theta}$ and $B_g^{\alpha, \theta}$
are \emph{asymmetric truncated Toeplitz operators} (ATTO) in $K_\theta$ (that is, general WH operators where
$P_1$ and $P_2$ are projections in $K_\theta$ and $A$ is a Toeplitz operator).

A natural motivation for studying ATTO comes from the study of finite interval convolution equations. Indeed, the close connection, via the Fourier transform, between TTO on a model space $K_\theta$ with $\theta(\xi)=e^{i\lambda\xi}\,,\,\lambda>0$ and finite interval convolution operators of the form
\[
W_\phi f (x)=\int_I f(t) \phi(x-t)\, dt\,\,,\,\,x\in I_1,
\]
where $I=I_1=[0,\lambda]$, has been highlighted, for instance, in \cite{BCFT}. It is clear that a similar relation exists between ATTO and finite interval convolution operators of the above form where $I_1\neq I$.

\bigskip

In Section~\ref{sec:2} we recall the definitions and basic properties of model spaces in an $H_p$ context,
while also introducing the notion of partial conjugation. Section~\ref{sec:3} analyses an isometric
isomorphism between $L_p$ spaces on the disc and half-plane, which restricts to $H^+_p$ and indeed
$\theta H^+_p$. For $p=2$ it has further properties which aid in the study of ATTO.
In Section~\ref{sec:4} ATTO are treated in some detail, and we solve
the question of uniqueness of symbol, via the characterization of the zero operator.
In Section~\ref{sec:5} we discuss the question which ATTO have finite rank.
Next, in Section~\ref{sec:6} it is shown that
ATTO are equivalent by extension to 
Toeplitz operators with triangular $2 \times 2$
matrix symbol.  This immediately enables one to obtain new results about ATTO (and even TTO)
from known results about standard Toeplitz operators.
Finally, Section~\ref{sec:7} discusses kernels of ATTO and the link with invariant subspaces.\\

%%%%%%%%%%%%%%%%%%%%           Section 2 - Model spaces, projections and Toeplitz operators         

\section{Model spaces, projections and Toeplitz operators}
\label{sec:2}

Recall that we write $L_p$ for $L_p(\RR)$, $H_p^+$  and $H_p^-$ for the Hardy
spaces of the upper and lower half-planes $\CC^+$ and $\CC^-$ (here $1 \le p \le \infty$) and we denote by $P^\pm$ the Riesz projections $P^+: L_p \to H_p^+$
and $P^-: L_p \to H_p^-$ for $1<p<\infty$.

For $\theta$ an inner function  (in $H_\infty^+$), let
$K_\theta=K_\theta^p$ denote the model space defined in \eqref{eq:defktheta},
where we omit the superscript $p$ unless it is necessary for clarity. If $\alpha$ and $\theta$ are inner functions, we say that $\alpha K_\theta$ is a \emph{shifted model space}. It is clear that $\alpha K_\theta \subset K_{\alpha \theta}$.

For any inner function $\theta$, we have the decomposition \eqref{eq:decomp}, and
  \begin{equation}\label{2}
   H_p^+=K_\theta \oplus \theta H_p^+\,,
\end{equation}
where the sum is orthogonal in the case $p=2$.
Let $P_\theta : L_p\rightarrow K_\theta$ be the projection from $L_p$ onto $K_\theta$ defined
by \eqref{eq:decomp}; we have
\beq \label{3}
P_\theta  = \theta P^- \overline\theta P^+ = P^+ \theta P^- \overline\theta I.
\eeq
Let moreover $Q_\theta$ be the operator defined in $L_p$, $1<p<\infty$, by
\beq \label{4}
Q_\theta:=P^+ - P_\theta,
\eeq
and let us use the same notation $P_\theta, Q_\theta$  for $P_{\theta\mid_{H_p^+}}, Q_{\theta\mid_{H_p^+}} $, respectively. 
For any $\phi\in H_p^+$, we define
\beq \label{1A}
\phi^\theta=P_\theta \,\phi.%\quad,\quad\phi_-^\theta=\overline\theta\,\phi^\theta\,.
\eeq

Now take $g \in L_\infty$. The \emph{Toeplitz operator} with symbol $g$ in $K_\theta$ is the operator
\[
T_g:H_p^+\rightarrow H_p^+\,\,,\,\,T_g=P^+gI_+
\]
where $I_+$ denotes the identity operator in $H_p^+$. This definition can be generalised to the vectorial case straightforwardly, for a matricial symbol $g\in L_\infty^{n\times n}$. \\

%%%%%%%%%%%%  Prop 2.1
%\begin{prop}\label{propI.1}
%The following relations hold:
%\beq \label{5}
%P^+=T_{\bar{\theta}}T_\theta
%\eeq
%\beq \label{6}
%Q_\theta=T_\theta T_{\bar{\theta}}
%\eeq
%\beq \label{7}
%P_\theta=P^+-T_\theta T_{\bar{\theta}}=T_{\bar{\theta}}T_\theta-T_\theta T_{\bar{\theta}}
%\eeq
%\beq \label{8}
%T_\theta T_{\bar{\theta}}T_\theta=T_\theta
%\eeq
%\beq \label{9}
%T_{\bar{\theta}} T_\theta T_{\bar{\theta}}=T_{\bar{\theta}}
%\eeq
%\beq \label{10}
%Q_\theta T_\theta= T_\theta= \theta P^+
%\eeq
%\beq \label{11}
%T_{\bar{\theta}}Q_\theta= T_{\bar{\theta}}= \bar{\theta}  Q_\theta
%\eeq
%\beq \label{12}
%T_{\bar{\theta}}P_\theta= 0
%\eeq
%\beq \label{13}
%P_\theta T_\theta= 0.
%\eeq
%\end{prop}

%%%%%%%%%%%%%%%%%%%%%

If $\alpha, \theta\in H_\infty^+$ are inner functions, we say that $\alpha \preceq \theta$ if and only if there exists an inner function $\tilde\theta$ such that $\theta=\alpha\tilde\theta$, and $\alpha \prec \theta$ if and only if $\tilde\theta$ is not constant. 
Of course
\beq \label{15}
\alpha \preceq \theta \Rightarrow\bar{\alpha} \theta \preceq \theta.
\eeq
We also have
\beq \label{14}
\alpha \preceq \theta \Leftrightarrow K_\alpha \subset K_\theta \Leftrightarrow \ker P_\theta \subset \ker P_\alpha.
\eeq

%%%%%%%%  Prop 2.2

%\begin{prop} \label{propI.2}
%If $\alpha \preceq \theta$, then the following relations hold:
%\beq \label{16}
%P_\alpha P_\theta =P_\theta P_\alpha= P_\alpha
%\eeq
%\beq \label{17}
%Q_\alpha Q_\theta =Q_\theta Q_\alpha= Q_\theta
%\eeq
%\beq \label{18}
% Q_\theta P_\alpha =P_\alpha Q_\theta = 0
%\eeq
%\beq \label{19}
% Q_\alpha P_\theta =P_\theta Q_\alpha
%\eeq
%\beq \label{20}
% T_\alpha T_{\bar{\theta}} =\alpha \bar{\theta} Q_\theta = T_{\alpha \bar{\theta}} Q_\theta
%\eeq
%\beq \label{21}
% T_{\bar{\theta}} T_\alpha = T_{\alpha \bar{\theta}}
%\eeq
%\beq \label{22}
% T_\theta T_{\bar{\alpha}} =\bar{\alpha} \theta Q_\alpha = T_{\bar{\alpha} \theta} Q_\alpha
%\eeq
%\beq \label{23}
% T_{\bar{\alpha}} T_\theta = T_{\bar{\alpha} \theta}= \bar{\alpha} \theta P^+\,.
%\eeq
%\end{prop}
%%%%%%%%%%%

As a consequence of this we can also define, for $\alpha\preceq \theta$, a projection in $L_p$ (or $H_p^+$) by
\beq \label{24}
  P_{\alpha, \theta}:=P_\theta - P_\alpha\,\,,
\eeq
and again we use the same notation for the operator defined by \eqref{24} in $L_p$ and its restriction to $H_p^+$. 
% From Proposition \ref{propI.2} 
We  easily see that
\beq \label{25}
  P_{\alpha, \theta}=Q_\alpha P_\theta = P_\theta\, Q_\alpha  
=\alpha P_{\bar{\alpha} \theta}\, \bar{\alpha} I
\eeq
and it follows from \eqref{25} that the image of $P_{\alpha, \theta}$ is the shifted model space
\beq
K_{\alpha,\theta}:=K_\theta\cap\alpha H_p^+=\alpha K_{\bar{\alpha} \theta}.
\eeq
Of course $K_{\alpha,\theta}=K_\theta$ if $\alpha$ is constant, and $K_{\alpha,\theta}=K_\theta\ominus K_\alpha$ if $p=2$. 

\vspace{0,3cm}

We introduce now a class of conjugate-linear operators in $H_p^+$ by generalising the notion of a conjugation in a complex Hilbert space $\cal{H}$  (i. e., an isometric conjugate-linear involution in $\cal{H}$).

%%Definition 2.3
\begin{defn}\label{def2.3}
 Let $X, Y$ be closed subspaces of $H_p^+$ such that $x\bot y$ for all $x\in X\cap H_2^+, y\in Y\cap H_2^+$, and let $A=X\oplus Y$. We say that a conjugate-linear operator in $H_p^+$, $\cal{C}$, is a {\rm partial conjugation} in $A$ if and only if $\cal{C}_{\mid_X}$ is an isometric involution on $X$ and $\cal{C}_{\mid_Y}$=$0$.

If $Y=\{0\}$ then $\cal{C}$ is a {\rm conjugation} in $A$.
\end{defn}

Let now $\cal C_\theta$ be the conjugate-linear operator defined in $H_p^+$, for each inner function $\theta$, by
\beq \label{c1}
{\cal C}_\theta (\phi_+)=\theta\, \overline{P_\theta\,\phi_+}\quad,\quad \phi_+\in H_p^+.
\eeq
It is easy to see that $({\cal C}_\theta)^2=P_\theta$, ${\cal C}_\theta$ maps $K_\theta$ onto $K_\theta$ isometrically, and ${\cal C}_\theta (\theta H_p^+)=\{0\}$. Thus ${\cal C}_\theta$ is a partial conjugation in $H_p^+$ and, analogously ${\cal C}_\alpha$ is a partial conjugation in $K_\theta$ if $\alpha \preceq \theta$. Of course ${\cal C}_\alpha$ is a conjugation in $K_\alpha$.

\vspace{0,3cm}

We will also use the following simple relations.
Let $r$ denote the function defined by
\beq\label{def:r}
r(\xi)=\frac{\xi-i}{\xi+i}
\eeq
for $\xi \in \CC$ and let $\phi_\pm\in H_p^\pm$. Then
\beq\label{A}
P^+r^{-1}\phi_+=r^{-1}\phi_+-2i\frac{\phi_+(i)}{\xi-i},
\eeq
\beq\label{B}
P^-r\phi_-=r\phi_-+2i\frac{\phi_-(-i)}{\xi+i}.
\eeq
Moreover, if $\theta$ is an inner function, taking into account that $\phi_+=\phi_+^\theta+\theta\tilde\phi_+$ with $\tilde\phi_+\in H_p^+$, we have

\beq\label{C}
P_\theta h_+\phi_+=P_\theta h_+\phi_+^\theta
\eeq
whenever $h_+$ is such that $h_+\phi_+\in L_p$ and $h_+Q_\theta\phi_+\in  \theta H_p^+$ (in particular, if $h_+\in H_\infty^+$), and
\beq\label{D}
Q_\theta h_-\phi_+^\theta=0\quad,\quad P_\theta h_-\phi_+^\theta=P^+ h_-\phi_+^\theta
\eeq
whenever $h_-$ is such that $h_-\phi_+^\theta\in L_p$ and $h_-\bar\theta\phi_+^\theta\in H_p^-$ (in particular, if $h_-\in H_\infty^-$).  As a consequence of \eqref{C} and \eqref{D}, we also have
\beq\label{E}
\alpha\preceq\theta\Rightarrow P_\theta h_-\phi_+^\alpha=P_\alpha h_-\phi_+^\alpha\quad,\quad P_\alpha \, h_+\phi_+^\theta=P_\alpha\,h_+\phi_+^\alpha.
\eeq

%% Section 3 - Equivalence between operators on the disc and half-plane
\section{Equivalence between operators on the disc and half-plane}
\label{sec:3}

We now recall the details of the isometric isomorphism between the Hardy spaces $H_p^+$
on the upper half-plane $\CC^+$ and $H_p(\DD)$ on the unit disc $\DD$. It will be seen that
this leads to an isometric bijective equivalence (i.e., an unitary equivalence in the case $p=2$) between
model spaces on the disc and half-plane; in the case $p=2$ this leads to a unitary equivalence between
(A)TTO on the disc and half-plane, enabling us to
give an immediate translation of our results to the disc context.

Our convention in this section is that lower case letters such as $f$   denote functions
on the disc, whereas capital letters denote functions on the half-plane. \\

Let $m: \DD \to \CC^+$ be the conformal bijection given by
\[
m(z)=i\left( \frac{1-z}{1+z} \right), \qquad m^{-1}(\xi)=\frac{i-\xi}{i+\xi},
\]
(other choices are possible)
and $V: H_p(\DD)\to H_p(\CC^+)$ the isometric isomorphism given by
\beq\label{eq:Viso}
(Vf)(\xi)= \frac{1}{\pi^{1/p}} \frac{1}{(i+\xi)^{2/p}} \, f(m^{-1}(\xi)), \qquad (f \in H_p(\DD)),
\eeq
(see, for example, 
\cite{Koosis,nik,jrp-hankel}).
The inverse mapping is given by
\[
(V^{-1}F)(z)=\pi^{1/p}\left(\frac{2i}{1+z}\right)^{2/p} F(m(z)), \qquad (F \in H_p(\CC^+)).
\]
Now for $n \in \ZZ$ the function $z^n$ is mapped by $V$ to the function $e_n$ given by
\[
e_n^{(p)}(\xi)=\frac{1}{\pi^{1/p}} \frac{(i-\xi)^n}{(i+\xi)^{n+2/p}}.
\]
The same formula (\ref {eq:Viso}) extends $V$ to an isometric mapping from $L_p(\TT)$
onto $L_p(\RR)$, and for $p = 2$ it also maps $\overline{H^p_0(\DD)}$ into $H_p^-$.

 Let $\theta$ be an inner function
in $H^\infty(\CC^+)$; then the function $\Theta:=\theta \circ m^{-1}$ is an inner function
in $H^\infty(\DD)$.
Now for $f=\theta g$ with $g \in H^p(\DD)$ we have
\[
(Vf)(\xi)=  \Theta(\xi) (Vg)(\xi),
\]
so $V$ takes $\theta H^p(\DD)$ onto $\Theta H^p(\CC^+)$.  

Letting $q$ be the conjugate index to $p$, we also have that $(V^*)^{-1}$ maps
$H_q(\DD)$ onto $H_q(\CC^+)$ and takes its subspace $K_\theta$ to $K_\Theta$.

The situation is better for $p=2$,  since $V$ 
is unitary, and it maps $K_\theta=H_2(\DD) \cap \theta\overline{H^2_0(\DD)}$
onto $K_\Theta= H_2^+ \cap \Theta H_2^-$; hence, the decomposition 
\[
L^2(\TT) =  \overline{H^2_0(\DD)} \oplus K_\theta \oplus \theta H^2(\DD)
\]
is mapped by $V$ term-wise into
\[
L^2(\RR) = H_2^- \oplus K_\Theta \oplus \Theta H_2^+.
\]
This situation does not hold for $p \ne 2$.\\

Suppose now that $p=2$ and $g \in L^\infty(\DD)$. We write $G:=g \circ m^{-1}$ and $\gA=\alpha \circ m^{-1}$.
Then, we
have the following commutative diagram, where $A^{\alpha,\theta}_g$ denotes an ATTO on the disc,
as defined analogously to \eqref{II.1}:

\beq\label{eq:commdiag}
\begin{matrix}
K_\theta & \stackrel{ A^{\alpha,\theta}_g}{\longrightarrow}& K_\alpha \\
V\downarrow && \downarrow V\\
K_{\Theta } & \stackrel{A^{\gA,\Theta}_{G}}{\longrightarrow} &K_{\gA}
\end{matrix}
\eeq
We see that this diagram commutes, 
since for $k \in K_\theta$ we have 
\begin{eqnarray*}
V(gk)(\xi) &=& \frac{1}{\pi^{1/2}} \frac{1}{(i+\xi)}  g(m^{-1}(\xi))k(m^{-1}(\xi)) \\
&=& G(\xi) (Vk)(\xi);
\end{eqnarray*}
now, since $P_\gA V=VP_\alpha$ 
we get 
\[
VP_\alpha (gk)=P_\gA V(gk) = P_{\gA}G(Vk),
\]
so we have the required unitary equivalence between ATTO on the disc and half-plane.

%%%%%%%%%%%%%
%%%%%% Section 4 - Asymmetric truncated Toeplitz operators

\section{Asymmetric truncated Toeplitz operators}
\label{sec:4}

Let $g\in L_\infty$ and let $\alpha,\theta\in H_\infty^+$ be inner functions. 
As in Section~\ref{sec:1}, we define the \emph{asymmetric truncated Toeplitz operators} (abbreviated to \emph{ATTO}) $A_g^{\alpha, \theta}$ and, for $\alpha\preceq\theta$, $B_g^{\alpha, \theta}$ as follows:  

\beq \label{II.2}
 A_g^{\alpha, \theta} = P_\alpha \,g\,P_\theta%_{\mid_{K_\theta}}
\eeq
 \beq \label{II.3}
 B_g^{\alpha, \theta} = P_{\alpha, \theta} \;g\,P_\theta%_{\mid_{K_\theta}}
\eeq
where $A_g^{\alpha, \theta}$ and $B_g^{\alpha, \theta}$ can be seen as operators in $H_p^+$, or operators in $K_\theta$ if $\alpha\preceq\theta$, or as operators from $K_\theta$ into $K_\alpha$ and $K_{\alpha, \theta}$, respectively. We will assume the latter unless stated otherwise. If $\alpha=\theta$ then $A_g^{\alpha,\theta}$ is the \emph{truncated Toeplitz operator} $A_g^\theta$.

It is easy to see that $A_g^\theta=A_g^{\alpha, \theta}+B_g^{\alpha, \theta}$ and that an ATTO of the form \eqref {II.3} can be expressed in terms of ATTO of type \eqref {II.2}, since we have  
\beq \label{II.5}
B_g^{\alpha, \theta} =P_{\alpha, \theta}T_{g\mid_{K_\theta}}  =\alpha P_{\bar{\alpha} \theta}\bar{\alpha}T_{g\mid_{K_\theta}}   
     =\alpha P_{\bar{\alpha} \theta}T_{\bar{\alpha}g\mid_{K_\theta}}  =\alpha A_{\bar{\alpha} g}^{\bar{\alpha} \theta, \theta}.
\eeq
We will therefore focus here on ATTO of type \eqref {II.2}. Moreover, considering that
\[
{(A_g^{\alpha, \theta})}^*=A_{\bar g}^{\theta,\alpha},
\]
we will assume in what follows that $\alpha\preceq\theta$.\\

\vspace{0,3cm}

We will use the following generalisation of the notion of a complex symmetric operator in a Hilbert space.

%%Definition 4.1
\begin{defn} Let $A$ be a closed subspace of $H_p^+$. An operator $T:A\rightarrow H_p^+$ is a {\rm complex partially symmetric operator} (respectively, a {\rm complex symmetric operator}) if and only if there exists a partial conjugation (respectively, a conjugation) in $A$, $\cal{C}$, such that $\cal{C} T\cal{C}=\tilde T$, where $\tilde T$ coincides with $T^*$ in $H_p^+\cap H_q^+,\, 1/p+1/q=1$. In this case we say that $T$ is \emph{$\cal{P}\cal{C}$-symmetric} (respectively, \emph{$\cal{C}$-symmetric}).
\end{defn}

%%Theorem 4.2
\begin{thm} \label{thm:3.2a}
If  $g\in L_\infty$, then
\[
{\cal C}_\alpha A_{g}^{\alpha, \theta}{\cal C}_\alpha=A_{\overline {g}}^\alpha\, .
\]
\end{thm}

\beginpf
Let $\phi_+\in K_\theta$. Then
\begin{eqnarray}\nonumber
{\cal C}_\alpha A_{g}^{\alpha, \theta}{\cal C}_\alpha\, \phi_+&=&\alpha \overline{A_{g}^{\alpha, \theta}{\cal C}_\alpha \phi_+}=\alpha \overline{P_\alpha g P_\theta \,\alpha \overline{P_\alpha\phi_+}}=\alpha\overline{P_\alpha g \alpha \overline{P_\alpha \phi_+}}\\
\nonumber &=&\alpha (\bar\alpha P^+ \alpha P^- \overline{g} \bar\alpha P_\alpha\phi_+)=P^+\alpha P^- \bar\alpha (P^++P^-)\bar g P_\alpha \phi_+\\
\nonumber &=&\alpha P^- \bar\alpha P^+\bar  g P_\alpha \phi_+=P_\alpha \bar g P_\alpha \phi_+,
\end{eqnarray}
for all $\phi_+\in H_p^+$.
\endpf

%%Corollary 4.3
\begin{cor}
For $g\in L_\infty$, $A_g^\theta$ is $\cal{C}_\theta$-symmetric in $K_\theta$ and we have
\beq\label{8bis}
{\cal C}_\theta A_{g}^\theta=A_{\overline {g}}^\theta\, {\cal C}_\theta.
\eeq
\end{cor}

Let us consider now the case of analytic symbols $g_+\in H_\infty^+$.

%%%%%%Theorem 4.4
\begin{thm} \label{prop:3.1a}
(i) If $g_+\in H_\infty^+$ and $\alpha, \theta$ are inner functions with $\alpha\preceq\theta$, then 
\[
A_{g_+}^{\alpha, \theta} \phi_+=A_{g_+}^\alpha \phi_+\,, \quad A_{\overline{g_+}}^{\theta,\alpha} \phi_+=A_{\overline{g_+}}^\alpha \,\phi_+
\]
for all $\phi_+\in H_p^+$.

(ii) If $\alpha\preceq\beta$ and $\beta\preceq\theta$, then $A_{g_+}^{\alpha, \beta} A_{f_+}^{\beta, \theta}=A_{g_+f_+}^{\alpha, \theta}$. 
\end{thm}

\beginpf
(i) follows from \eqref{E}.

(ii)$A_{g_+}^{\alpha, \beta} A_{f_+}^{\beta, \theta}=P_\alpha g_+ P_\beta f_+ P_\theta=P_\alpha g_+(P^+- Q_\beta) f_+ P_\theta=P_\alpha g_+ f_+ P_\theta=A_{g_+f_+}^{\alpha, \theta}$.
\endpf

As an immediate consequence we have, for $g_+\in H_\infty^+, n\in \mathbb N$, 
\beq \label{p1}
{(A_{g_+}^\theta)}^n=A_{g_+^n}^\theta.
\eeq

From Theorems \ref{thm:3.2a} and \ref{prop:3.1a} we also have the following.

%%%%%%Theorem 4.5
\begin{thm} \label{prop:3.1abis}
If $g_+\in H_\infty^+$, then $A_{g_+}^{\alpha, \theta}$ and $A_{\overline {g_+}}^{\theta,\alpha}$ are $\cal{P}\cal{C}_\alpha$-symmetric and 
\[
{\cal C}_\alpha A_{g_+}^{\alpha, \theta}=A_{\overline {g_+}}^\alpha\,{\cal C}_\alpha =A_{\overline{g_+}}^{\theta,\alpha}\,{\cal C}_\alpha .
\]
\end{thm}

\beginpf
By Theorem \ref{thm:3.2a} we have ${\cal C}_\alpha A_{g_+}^{\alpha, \theta}=A_{\overline {g_+}}^\alpha\,{\cal C}_\alpha$ and, by Theorem \ref{prop:3.1a} (i), $A_{\overline{g_+}}^{\theta,\alpha}=A_{\overline {g_+}}^\alpha$.
\endpf

Let us now consider the functions $k_w^\theta$ and $\tilde k_w^\theta$ defined, for each $w\in \mathbb C^+$, by
\beq\label{k1}
k_w^\theta (\xi):= \frac{1-\overline{\theta (w)}\,\theta (\xi)}{\xi - \overline w}\,,
\eeq
\beq\label{k2}
\tilde k_w^\theta (\xi):= \frac{\theta (\xi)-\theta (w)}{\xi - w}\,,
\eeq
which will play an important role in this section.
We have $k_w^\theta, \tilde k_w^\theta \in K_\theta$, with
\beq\label{eq:kk1311}
k_w^\theta=P_\theta \,\frac{1}{\xi-\overline w}\quad,\quad \tilde k_w^\theta =P_\theta\, \frac{\theta}{\xi-w}={\cal C}_\theta k_w^\theta.
\eeq
If $\alpha\preceq\theta$, the functions $k_w^\alpha, \tilde k_w^\alpha$ are related to $k_w^\theta, \tilde k_w^\theta$, respectively, by
\beq\label{k10}
P_\alpha k_w^\theta= k_w^\alpha\quad,\quad P_\alpha \tilde k_w^\theta=(\bar\alpha \theta)(w) \tilde k_w^\alpha\,.
\eeq

%%%%%%Theorem 4.6
\begin{thm} \label{thm:3.3a}
$k_i^\theta$ is a cyclic vector for $A_r^\theta$ and $\tilde k_i^\theta$ is a cyclic vector for $A_{r^{-1}}^\theta$.
\end{thm}

\beginpf
By \eqref{p1} and \eqref{C},
\[
(A_r^\theta)^n\,k_i^\theta=A_{r^n}^\theta\,k_i^\theta=P_\theta \,r^n P_\theta\, \frac{1}{\xi+i}=P_\theta \,(r^n\frac{1}{\xi+i}),
\]
so $\{(A_r^\theta)^n\,k_i^\theta\,:\,n\in\mathbb N\}$ is dense in $K_\theta$. On the other hand, 
since $T_\theta T_{\bar{\theta}}T_\theta=T_\theta$, we have
% by \eqref {8},
\[
(A_{r^{-1}}^\theta)^n\,\tilde k_i^\theta=A_{r^{-n}}^\theta\,{\cal C}_\theta\, k_i^\theta={\cal C}_\theta\,A_{r^{n}}^\theta\,k_i^\theta
\]
and, since ${\cal C}_\theta$ is an isometry in $K_\theta$, it follows that $\tilde k_i^\theta$ is a cyclic vector for $A_{r^{-1}}^\theta$.
\endpf

%%%%%%Theorem 4.7
\begin{thm} \label{thm:3.4a}
The operators $P_\alpha - A_r^{\alpha, \theta}\,A_{r^{-1}}^{\theta, \alpha}$ and $P_\alpha - A_{r^{-1}}^{\theta,\alpha}\,A_r^{\alpha, \theta}\,$ on $H+p^+$ are rank-one operators, with range equal to $\spam\{k_i^\alpha\}$and $\spam\{\tilde k_i^\alpha\}$, respectively, and we have
\beq\label{kVI.b}
(P_\alpha - A_r^{\alpha, \theta}\,A_{r^{-1}}^{\theta, \alpha})\,\phi_+=2i\phi_+^\alpha (i) \,k_i^\alpha\,,
\eeq
\beq\label{kVI.d}
(P_\alpha - A_{r^{-1}}^{\theta,\alpha}\,A_r^{\alpha, \theta})\,\phi_+=-2i\phi_-^\alpha (-i) \,\tilde k_i^\alpha\,,
\eeq
where $\phi_-^\alpha=\bar\alpha \phi_+^\alpha=\overline{{\cal C}_\alpha \phi_+^\alpha}$ .
\end{thm}

\beginpf
\begin{eqnarray}\nonumber
A_r^{\alpha, \theta}\,A_{r^{-1}}^{\theta, \alpha}\,\phi_+&=&P_\alpha\,rP_\theta\,r^{-1}\,\phi_+^\alpha=P_\alpha\,rP^+r^{-1}\,\phi_+^\alpha\\
\nonumber &=&\phi_+^\alpha-2i\phi_+^\alpha(i)\,P_\alpha \frac{1}{\xi+i}=\phi_+^\alpha-2i\phi_+^\alpha(i)\,k_i^\alpha,
\end{eqnarray}
where we used \eqref{A}, and \eqref{kVI.b} follows from this equality. 
On the other hand, by Theorem \ref{prop:3.1a}, Theorem \ref{prop:3.1abis}, \eqref{eq:kk1311} and \eqref{kVI.b},
% \eqref{E}, \eqref{8} and \eqref{kVI.b},
\begin{eqnarray}\nonumber
A_{r^{-1}}^{\theta, \alpha}\,A_r^{\alpha, \theta}\,\phi_+&=&A_{r^{-1}}^\alpha\,({\cal C}_\alpha)^2\,A_r^\alpha\,\phi_+={\cal C}_\alpha \,A_r^\alpha\,A_{r^{-1}}^\alpha\,{\cal C}_\alpha\,\phi_+\\
\nonumber &=&{\cal C}_\alpha A_r^{\alpha, \theta}\,A_{r^{-1}}^{\theta, \alpha}{\cal C}_\alpha\,\phi_+=-2i\,\alpha\,{\overline{({\cal C}_\alpha\,\phi_+)}_{(i)}}\,\overline {k_i^\alpha}=-2i{\phi_-^\alpha}(-i) \,\tilde k_i^\alpha .
\end{eqnarray}
\endpf

In particular, for $\alpha=\theta$, we have the defect operators (\cite{sarason07}) $I_{K_\theta}- A_r^\theta\,A_{r^{-1}}^\theta$ and $I_{K_\theta}- A_{r^{-1}}^\theta\,A_r^\theta$ in $K_\theta$, where $ I_{K_\theta}$ denotes the identity operator in $K_\theta$, with
\beq\label{VI.e}
(I_{K_\theta}- A_r^\theta\,A_{r^{-1}}^\theta)\,\phi_+^\theta=2i\, \phi_+^\theta(i)\,k_i^\theta
\eeq
\beq\label{VI.e1}
(I_{K_\theta}- A_{r^{-1}}^\theta\,A_r^\theta)\,\phi_+^\theta=-2i\, \phi_-^\theta(-i)\,\tilde k_i^\theta .
\eeq

\vspace{0,2cm}

Next we address the question when an ATTO is zero, which is equivalent to obtaining conditions for two ATTO to be equal. 
For this purpose, it will be useful to note that a symbol $g\in L_\infty$ admits the following decompositions:
\beq\label{k13}
g=G_++G_-\;, {\rm with} \; \;G_\pm=(\xi +i) P^\pm\frac{g}{\xi+i}\;,
\eeq
\beq\label{k14}
g=g_++g_-\;, {\rm with} \; \;g_\pm=(\xi -i) P^\pm\frac{g}{\xi-i}\;,
\eeq
\beq\label{k15}
g=\gamma_++\gamma_-+C\;, {\rm with} \; \;\gamma_\pm=(\xi \pm i) P^\pm\frac{g}{\xi\pm i}\;,\;C\in \CC.
\eeq

The third decomposition can easily be related to any of the other two; for instance,
\[
G_+ =\gamma_+\,,\,G_-=\gamma_-+C \;{\rm with}\; C=-2iP^-(\frac{g}{\xi-i})\,(-i).
\]

It is clear that an ATTO does not have a unique symbol, since we can have $A_g^{\alpha, \theta}=0$ with $g\neq 0$. In fact, using the previous results and defining ${\cal H}_p^\pm:=\lambda_\pm H_p^\pm$ where $\lambda_\pm (\xi)=\xi\pm i$, we have the following.

%%%%%Theorem 4.8
\begin{thm} \label{thm:3.1a}
$A_g^{\alpha, \theta}=0$ if and only if $g=\bar\theta\tilde g_-+\alpha\tilde g_+$ with $\tilde g_\pm\in {\cal H}_p^\pm$.
\end{thm}

\beginpf
First we prove that $A_g^{\alpha, \theta}=0$ if $g=\bar\theta\tilde g_-+\alpha\tilde g_+$.
For $z_+\in\mathbb C^+$, let $k_{z_+}^\theta:=\frac{1-\overline{\theta(z_+)}\theta}{\xi -\overline{z_+}}=P_\theta(\frac{1}{\xi -\overline{z_+}})$; then
\begin{eqnarray}
A_g^{\alpha, \theta} k_{z_+}^\theta&=&P_\alpha \,[g\, \frac{1-\overline{\theta(z_+)}\theta}{\xi -\overline{z_+}}]=P_\alpha \,[(\bar\theta\tilde g_-+\alpha\tilde g_+)\, \frac{1-\overline{\theta(z_+)}\theta}{\xi -\overline{z_+}}]\nonumber\\
&=&P_\alpha \,[\tilde g_-\, \frac{\bar\theta-\overline{\theta(z_+)}}{\xi -\overline{z_+}}]+P_\alpha \,[\alpha\tilde g_+\, \frac{1-\overline{\theta(z_+)}\theta}{\xi -\overline{z_+}}]=0\nonumber
\end{eqnarray}
since $\tilde g_-\, \frac{\bar\theta-\overline{\theta(z_+)}}{\xi -\overline{z_+}}\in H_p^-$ and $\alpha\tilde g_+\, \frac{1-\overline{\theta(z_+)}\theta}{\xi -\overline{z_+}}\in\alpha H_p^+$.
The converse will be proved in several steps. Assuming that $A_g^{\alpha, \theta}=0$, we show that

\vspace{0,1cm}
(i) $A_{G_+}^{\alpha,\theta}A_r^{\alpha,\theta}A_{r^{-1}}^{\theta,\alpha}\,k_i^\alpha=A_r^{\alpha,\theta}A_{r^{-1}}^{\theta,\alpha}A_{G_+}^{\alpha,\theta}\,k_i^\alpha$,

\vspace{0,1cm}

(ii) $\gamma_+=\alpha f_++C_1$\quad  for some $f_+\in{\cal H}_p^+$ and some $C_1\in\CC$,

\vspace{0,1cm}

(iii) $\gamma_-=\bar\theta f_-+C_2$\quad  for some $f_-\in{\cal H}_p^-$ and some $C_2\in\CC$,

\vspace{0,1cm}

(iv) $C_1+C_2+C=0$, where $C$ is the constant in \eqref{k15},

\vspace{0,1cm}

so that $g=\alpha f_++\bar\theta f_-$ with $f_\pm\in{\cal H}_p^\pm$.

\vspace{0,2cm}

(i) Let $G_\pm$ be defined as in \eqref{k13}. 
We have, from \eqref{kVI.b},
\[
A_{G_+}^{\alpha,\theta}A_r^{\alpha,\theta}A_{r^{-1}}^{\theta,\alpha}\,k_i^\alpha=(1-2ik_i^\alpha(i)) P_\alpha G_+\,k_i^\alpha.
\]
Now, if $A_g^{\alpha, \theta}=0$ then $A_{G_++G_-}^{\alpha, \theta}=0$ and 
\beq\label{k16}
A_{G_+}^{\alpha, \theta}\,\phi_+=-A_{G_-}^{\alpha, \theta}\,\phi_+ 
\eeq
for all $\phi_+$ such that $G_\pm\phi_+^\theta\in H_p^+$ (where we define $A_{G_\pm}^{\alpha, \theta}\,\phi_+=P_\alpha\,G_\pm\phi_+^\theta)$. Also note that
\beq\label{k16a}
P_\alpha G_-k_i^\alpha=P^+ G_-k_i^\alpha.
\eeq
Using \eqref{k16}, \eqref{k16a}, \eqref{D}, \eqref{E}, and taking into account that
\[
P_\theta r^{-1}k_i^\alpha=P_\alpha r^{-1}k_i^\alpha=P^+ r^{-1}k_i^\alpha=r^{-1}k_i^\alpha -2i\,\frac{k_i^\alpha (i)}{\xi-i},
\]
we have
\begin{eqnarray}\nonumber
A_r^{\alpha,\theta}A_{r^{-1}}^{\theta,\alpha}A_{G_+}^{\alpha,\theta}\,k_i^\alpha&=& -A_r^{\alpha,\theta}A_{r^{-1}}^{\theta,\alpha}A_{G_-}^{\alpha,\theta}\,k_i^\alpha=-A_r^{\alpha,\theta}(P_\alpha r^{-1}G_-\,k_i^\alpha)\\
\nonumber &=&-A_r^{\alpha,\theta}[P_\alpha G_-(r^{-1} k_i^\alpha- 2i\,\frac{k_i^\alpha (i)}{\xi-i})]=-A_r^\alpha A_{G_-}^{\alpha,\theta}r^{-1} k_i^\alpha\\
\nonumber &=&A_r^\alpha A_{G_+}^{\alpha,\theta}r^{-1} k_i^\alpha=P_\alpha\,rP_\alpha G_+P_\theta \,r^{-1} k_i^\alpha=P_\alpha rG_+P_\theta r^{-1} k_i^\alpha\\
\nonumber&=&P_\alpha rG_+P_\alpha r^{-1} k_i^\alpha=P_\alpha G_+ k_i^\alpha-2ik_i^\alpha(i)\,P_\alpha G_+ P_\alpha \frac{1}{\xi+i}\\
\nonumber&=&(1-2ik_i^\alpha(i))P_\alpha G_+k_i^\alpha.
\end{eqnarray}
Thus, $A_{G_+}^{\alpha,\theta}A_r^{\alpha,\theta}A_{r^{-1}}^{\theta,\alpha}\,k_i^\alpha=A_r^{\alpha,\theta}A_{r^{-1}}^{\theta,\alpha}A_{G_+}^{\alpha,\theta}\,k_i^\alpha$.

\vspace{0,2cm}

(ii) From (i) we get
\[
(A_{G_+}^{\alpha,\theta}-A_{G_+}^{\alpha,\theta}A_r^{\alpha,\theta}A_{r^{-1}}^{\theta,\alpha})\,k_i^\alpha=(A_{G_+}^{\alpha,\theta}-A_r^{\alpha,\theta}A_{r^{-1}}^{\theta,\alpha}A_{G_+}^{\alpha,\theta})\,k_i^\alpha
\]
and thus
\[
A_{G_+}^{\alpha,\theta}(P_\alpha -A_r^{\alpha,\theta}A_{r^{-1}}^{\theta,\alpha})\,k_i^\alpha=(P_\alpha -A_r^{\alpha,\theta}A_{r^{-1}}^{\theta,\alpha})A_{G_+}^{\alpha,\theta}\,k_i^\alpha
\]
which, by Theorem \ref {thm:3.4a}, is equivalent to
\[
A_{G_+}^{\alpha,\theta}2i\,k_i^\alpha(i)\,k_i^\alpha=2i(A_{G_+}^{\alpha,\theta}\,k_i^\alpha)(i)\,k_i^\alpha.
\]
Therefore,
\[
A_{G_+}^{\alpha,\theta}\,k_i^\alpha=C_1 k_i^\alpha, \quad {\rm where}\;  C_1\in \CC \setminus \{0\}
\]
and 
\[
A_{G_+}^{\alpha,\theta}\,k_i^\alpha=C_1 k_i^\alpha \Leftrightarrow P_\alpha (G_+-C_1) k_i^\alpha=0\Leftrightarrow P_\alpha \frac{G_+-C_1}{\xi+i}=0\Leftrightarrow \frac{G_+-C_1}{\xi+i}\in \alpha H_p^+.
\]
Since $G_+ =\gamma_+$, we have $\gamma_+=\alpha f_++C_1$\quad  with $f_+\in{\cal H}_p^+$ and $C_1\in\CC$.

\vspace{0,2cm}

(iii) Since $\bar g=(\bar g)_++(\bar g)_-$, where
\[
(\bar g)_\pm=(\xi+i) P^\pm \frac{\bar g}{\xi+i},
\]
so that $(\bar g)_+=\overline{\gamma_-}$, to study the condition on $\gamma_-$ we use the equivalence $A_g^{\alpha,\theta}=0\Leftrightarrow A_{\bar g}^{\theta,\alpha}=0\Leftrightarrow P_\theta\,\bar gP_\alpha=0$, where the equality on the right-hand side means that

\beq\label{k17}
P_\alpha\,\bar gP_\alpha=0\quad \wedge \quad P_{\alpha,\theta}\,\bar gP_\alpha=0.
\eeq

From the first equality in \eqref{k17} anf from (ii) we conclude that, for some constant $C_2\in\CC$,
\beq\label{k19}
\frac{(\bar g)_+-C_2}{\xi+i}\in \alpha H_p^+.
\eeq
On the other hand we have, from the second equality in \eqref{k17},
\beq\label{k20}
P_{\alpha,\theta}\,(\bar g)_+\, k_i^\alpha =-P_{\alpha,\theta}(\bar g)_-\, k_i^\alpha=-\alpha P_{\bar\alpha \theta}\bar\alpha(I-P^-)(\bar g)_-\,k_i^\alpha=0. 
\eeq
Since we also have $P_{\alpha,\theta}\,C_2 k_i^\alpha=0$, taking this and \eqref{k20} into account we get
\[
0=P_{\alpha,\theta}((\bar g)_+-C_2) k_i^\alpha=P_{\alpha,\theta}(\frac{(\bar g)_+-C_2}{\xi+i}(1-\overline{\alpha(i)} \alpha))
\]
which, by \eqref{k19}, implies that
\[
0=P_\theta(f+(1-\overline{\alpha(i)} \alpha))\quad {\rm with}\quad f_+=\frac{(\bar g)_+-C_2}{\xi+i}.
\]
Now,
\[
P_\theta [f_+(1-\ol {\alpha(i)} \alpha)]=0 \Rightarrow P_\theta f_+=0
\]
because $P_\theta [f_+(1-\ol {\alpha(i)} \alpha)]=0$ implies that $f_+(1-\ol {\alpha(i)}\alpha) = \theta \tilde f_+$, with $\tilde f_+\in H_p^+$ and, if $\tilde f_+î, \tilde f_+^o$ are the inner and outer factors of $\tilde f_+$, respectively, that is equivalent to having $f_+^i f_+^o(1-\ol {\alpha(i)}\alpha) = \theta \tilde f_+$.
Since $1-\ol {\alpha(i)}\alpha $ is an outer function in $H_\infty^+$, we conclude that $\theta$ divides $f_+^i$ and thus $P_\theta f_+=0$.
Thus $f_+\in\theta H_p^+$ and we conclude that $\gamma_-=\overline{(\bar g)_+}=\bar\theta f_-+C_2$ with $f_-\in {\cal H}_p^-$.

\vspace{0,2cm}

(iv) It follows from (ii), (iii) and \eqref{k15} that $g=\alpha f_++\bar\theta f_-+B$ where $B$ is a constant. Since $A_g^{\alpha,\theta}=0$, it follows from the first part of the proof that we must then have $A_B^{\alpha,\theta}=0$, which implies that $B=0$.
\endpf

For $p=2$ we may use the unitary equivalence derived earlier to obtain a generalisation of
Sarason's result for TTO in \cite{sarason07}, which, it seems, cannot be proved
directly using his techniques. It seems natural to conjecture that an analogous result holds in the disc for all
$1<p<\infty$, although no direct translation of the half-plane result seems to be possible for $p \ne 2$.

%% Corollary 4.9
\begin{cor}
In the case of $p=2$ and for Hardy spaces on $\DD$,
the asymmetric truncated Toeplitz operator  
$A_g^{\alpha, \theta}$ is zero
if and only if $g \in \alpha H^2(\DD) + \overline\theta \overline{H^2(\DD)}$.
\end{cor}
\beginpf
Note that
$g \in \alpha H^2(\DD) $ if and only if $g \circ m^{-1} \in (\alpha \circ m^{-1})\lambda_+H_2^+$
and $g \in \overline\theta \overline{H^2(\DD)}$ if and only if
$g \circ m^{-1} \in (\theta \circ m^{-1})\lambda_- H_2^-$.
Now the result follows directly from Theorem \ref{thm:3.1a}
using  the equivalence given in (\ref{eq:commdiag}).
\endpf

%%%%%%Section 5 - Finite rank asymmetric truncated Toeplitz operators
\section{Finite rank asymmetric truncated Toeplitz operators}
\label{sec:5}

In this section we assume again that $\alpha,\theta$ are inner functions with $\alpha \preceq \theta $. 
It is clear from any of the decompositions (\ref{k13})--(\ref{k15}) of $g \in L_\infty$ that we can represent
$g$ in the form
\begin{equation}\label{eq:F1}
g=a_+\overline\theta + a_- \alpha
\end{equation}
with $a_{\pm} \in {\mathcal H}_p^{\pm}$. If $a_{\pm} \in {\mathbb C}$, then by Theorem~\ref{thm:3.1a} we have
$A_g^{\alpha,\theta}=0$. It now seems natural to consider symbols of the form
\begin{equation}\label{eq:F2}
g=\frac{\alpha}{\xi-z_+} \quad \hbox{and} \quad g=\frac{\overline\theta}{\xi-\overline{z_+}}, \qquad (z_+ \in {\mathbb C}^+)
\end{equation}
as being the simplest corresponding to a non-zero ATTO $A_g^{\alpha,\theta}$.

Some other symbols seem equally simple. Let $\theta$ have a non-tangential limit $\theta(\xi_0)$ at $\xi_0 \in {\mathbb R}$
and suppose, in addition, that the functions
\begin{equation}\label{eq:F3}
\frac{\alpha(\xi)-\alpha(\xi_0)}{\xi-\xi_0} \quad \hbox{and} \quad \frac{\theta(\xi)-\theta(\xi_0)}{\xi-\xi_0} \quad \hbox{lie in} \quad L_\infty,
\end{equation}
in which case the functions in (\ref{eq:F3}) lie in $K_\alpha$ and $K_\theta$ respectively, and hence in $H_p^+$. We
can then consider bounded symbols of the form (\ref{eq:F1}) with
\[
a_-= \frac{\overline\theta-\overline{\theta(\xi_0)}}{\xi-\xi_0}, \quad a_+=\frac{\alpha(\xi_0)-\alpha}{\xi-\xi_0},
\]
i.e.,
\begin{equation}\label{eq:F4}
g=\frac{\alpha(\xi_0)\overline\theta-\overline{\theta(\xi_0)}\alpha}{\xi-\xi_0}.
\end{equation}

Analogously, if $\theta$ admits a non-tangential limit $\theta(\infty)$ at $\infty$, i.e., the 
inner function
$\theta\left(i\dfrac{1+z}{1-z}\right)$ in the unit disc has a non-tangential limit $\theta(\infty)$ at $1$, and in addition
the functions
\begin{equation}\label{eq:F5}
\xi[\alpha(\xi)-\alpha(\infty)] \quad \hbox{and} \quad \xi[\theta(\xi)-\theta(\infty)] \quad \hbox{lie in} \quad L_\infty,
\end{equation}
then we can consider bounded symbols of the form
\begin{equation}\label{eq:F6}
g=\xi[\alpha(\infty)\overline\theta-\overline{\theta(\infty)}\alpha].
\end{equation}
We remark that, if (\ref{eq:F5}) holds, then
\[
\tilde k^\alpha_\infty:=\alpha-\alpha(\infty) \in K_\alpha \quad \hbox{and}
\quad 
\tilde k^\theta_\infty:=\theta-\theta(\infty) \in K_\theta.
\]

%%Theorem 5.1
\begin{thm}\label {F1}
The asymmetric truncated Toeplitz operators $A_g^{\alpha,\theta}$ with $g$ of the form
(\ref{eq:F2}),
(\ref{eq:F4})
and 
(\ref{eq:F6}), are rank-one operators.
\end{thm}

\beginpf
Suppose that $g=\dfrac{\alpha}{\xi-z_+}$ with $z_+ \in {\mathbb C}^+$. Then for any $w \in {\mathbb C}^+$ and $k_w^\theta$ 
given by (\ref{k1}), we have
\begin{eqnarray*}
A_g^{\alpha,\theta}k_w^\theta &=& \alpha P^- \overline\alpha P^+ \frac{\alpha}{\xi-z_+} k_w^\theta
= \alpha P^- \overline\alpha \left( \frac{\alpha k_w^\theta-\alpha(z_+)k_w^\theta(z_+)}{\xi-z_+}\right)\\
&=& k_w^\theta(z_+) \frac{\alpha-\alpha(z_+)}{\xi-z_+}= k_w^\theta(z_+) \tilde k^\alpha_{z_+},
\end{eqnarray*}
where  $\tilde k^\alpha_{z_+}$ is defined in (\ref{k2}).

Analogously, if $g=\dfrac{\overline\theta}{\xi-\overline{z_+}}$ with $z_+ \in {\mathbb C}^+$, then
\[
A^{\alpha,\theta}_g k^\theta_w = - (\overline\theta k_w^\theta)(z_+) k^\theta_{z_+}
\]
for all $w \in \CC^+$.\\

Suppose now that $g$ takes the form (\ref{eq:F4}). Then, taking into account the fact that, for all $w \in \CC^+$
\[
\frac{k_w^\theta-k_w^\theta(\xi_0)}{\xi-\xi_0}
= \left( C_1+C_2 \frac{\theta-\theta(\xi_0)}{\xi-\xi_0} \right)\frac{1}{\xi-\overline w} \in H_p^+,
\]
where 
\[
C_1 = \frac{1+\overline{\theta(w)}\theta(\xi_0)}{\xi_0-\overline w} \quad \hbox{and} \quad
C_2 = -\frac{\xi_0 \overline{\theta(w)}}{\xi_0-\overline w},
\]
and
\[
\frac{\overline\theta k_w^\theta-(\overline\theta k_w^\theta)(\xi_0)}{\xi-\xi_0}
=
\left(\widetilde C_1+\widetilde C_2 \frac{\overline\theta-\overline{\theta(\xi_0)}}{\xi-\xi_0} \right)\frac{1}{\xi-\overline w} \in H_p^-,
\]
where
\[
\widetilde C_1 = \frac{\overline{\theta(w)}-\overline{\theta(\xi_0)}}{\xi_0-w} \quad \hbox{and} \quad
\widetilde C_2 = \frac{\xi_0-\overline w}{\xi_0-w},
\]
we have
\begin{eqnarray*}
A_g^{\alpha,\theta} k_w^\theta &=& P_\alpha \frac{\alpha(\xi_0)\overline\theta-\overline{\theta(\xi_0)}\alpha}{\xi-\xi_0}k_w^\theta \\
&=& P_\alpha 
\bigg[ \alpha(\xi_0)\underbrace{ \frac{\overline\theta k_w^\theta-(\overline\theta k_w^\theta)(\xi_0)}{\xi-\xi_0}}_{\in H_p^-}
-\overline\theta(\xi_0)\alpha \underbrace{\frac{k_w^\theta-k_w^\theta(\xi_0)}{\xi-\xi_0}}_{\in H_p^+}
-(\overline\theta k_w^\theta)(\xi_0) \underbrace{\frac{\alpha-\alpha(\xi_0)}{\xi-\xi_0}}_{\in K_\alpha}
\bigg] \\
&=& -(\overline\theta k_w^\theta)(\xi_0) \tilde k^\alpha_{\xi_0},
\end{eqnarray*}
where 
$\tilde k^\alpha_{\xi_0}:= \dfrac{\alpha-\alpha(\xi_0)}{\xi-\xi_0}$.\\

Let now $g$ take the form (\ref{eq:F6}). Then, for all $w \in \CC^+$, we have
\begin{eqnarray*}
A_g^{\alpha,\theta} k_w^\theta &=& P_\alpha \left[ \xi[\alpha(\infty)\overline\theta-\overline{\theta(\infty)}\alpha] \frac{1-\overline{\theta(w)}\theta}{\xi-\overline w} \right] \\
&=& P_\alpha
\bigg[
\alpha(\infty)\frac{\xi(\overline\theta-\overline{\theta(w)})}{\xi-\overline w}
+ \frac{\overline{\theta(w)}-\overline{\theta(\infty)}}{\xi-\overline w} \xi[\alpha-\alpha(\infty)]
+ \overline{\theta(\infty)}\overline{\theta(w)} \alpha \underbrace{\frac{\xi(\theta-\theta(\infty))}{\xi-\overline w}}_{\in H_p^+}
\bigg]\\
&=& \alpha P^- \overline\alpha P^+ \left[ \alpha(\infty) \frac{\xi(\overline\theta-\overline{\theta(\infty)})}{\xi-\overline w}\right]
+ \alpha P^-\overline \alpha \frac{\xi [\alpha-\alpha(\infty)]}{\xi-\overline w}(\overline{\theta(w)}-\overline{\theta(\infty)}) \\
&=& \alpha(\infty) \alpha P^- \overline\alpha \left( \overline w \frac{\overline{\theta(w)}-\overline{\theta(\infty)}}{\xi-\overline w}\right)\\
&&+ (\overline{\theta(w)}-\overline{\theta(\infty)}) \alpha(\infty) \alpha \left[
\frac{\xi(\overline{\alpha(\infty)}-\overline\alpha)-\overline w(\overline{\alpha(\infty)}-\overline{\alpha(w)})}{\xi-\overline w}
\right]\\
&=& \alpha(\infty)(\overline{\theta(w)}-\overline{\theta(\infty)}) \left[\overline w \frac{1-\overline{\alpha(w)}\alpha}{\xi-\overline w}
+
\frac{\alpha\xi(\overline{\alpha(\infty)}-\overline\alpha)-\alpha \overline w \overline{\alpha(\infty)}+\alpha \overline w \overline{\alpha(w)}}{\xi-\overline w}
\right] \\
&=& \alpha(\infty)(\overline{\theta(w)}-\overline{\theta(\infty)})(-1+\alpha \overline{\alpha(\infty)})\\
 &=& (\overline{\theta(w)}-\overline{\theta(\infty)})(\alpha-\alpha(\infty)) = (\overline{\theta(w)}-\overline{\theta(\infty)})k_w^\alpha.
\end{eqnarray*}
Since the span of $\{ k^\theta_w: w \in \CC^+\}$ is dense in $K_\theta$, we have proved the result.

\endpf

One can show analogously that if

(i) $g=\dfrac{\alpha}{(\xi-z_+)^n}$ or $g=\dfrac{\overline\theta}{(\xi-\overline{z_+})^n}$, with $n \in {\mathbb N}$, or

(ii) $\theta,\theta',\ldots,\theta^{(n-1)}$ have non-tangential limits
at $\xi_0 \in {\mathbb R}$, while
the functions $a_+$ and $a_-$ given by
\[
a_+(\xi)= \frac{ \alpha(\xi) - \sum_{j=0}^{n-1} \alpha^{(j)}(\xi_0)(\xi-\xi_0)^j/j!}{(\xi-\xi_0)^n}
\]
and
\[
\overline{a_-(\xi)}=\frac{ \theta(\xi) - \sum_{j=0}^{n-1} \theta^{(j)}(\xi_0)(\xi-\xi_0)^j/j!}{(\xi-\xi_0)^n}
\]
lie in $L_\infty$, and $g=a_+\overline\theta + a_- \alpha$,
or

(iii) $\theta,\theta',\ldots,\theta^{(n-1)}$ have non-tangential limits
at $\infty$, while
the functions $a_+$ and $a_-$ satisfying
\[
a_+(\xi)= \xi^n[a(\xi)-\sum_{j=0}^{n-1} a^{(j)}(\infty) \xi^{-j}/j!]
\]
and
\[
\overline{a_-(\xi)} = \xi^n [\theta(\xi)-\sum_{j=0}^{n-1} \theta^{(j)}(\infty) \xi^{-j}/j!]
\]
lie in $L_\infty$, and $g=a_+ \overline\theta + a_- \alpha$,
 then
$A^{\alpha,\theta}_g$ is a finite-rank operator. 

\vspace{0,3cm}

Whether every rank-one ATTO with symbol in $L_\infty$ is of the form considered in Theorem \ref{F1}, or every finite-rank ATTO with symbol in $L_\infty$ is a linear
combination of those given above is an open question.
For the space $H_2(\DD)$, the answer is affirmative, as follows from the work of Sarason \cite{sarason07}
and Bessonov~\cite{bessonov}.

%%%%%%%%%%%%%%%       Section 6-    Equivalence after extension of ATTO and $T$-operators with a  triangular matrix symbol

\section{Equivalence after extension of ATTO and $T$-operators with a triangular matrix symbol}
\label{sec:6}

In this section we show that asymmetric truncated Toeplitz operators are equivalent after extension to Toeplitz operators with triangular symbols of a certain form.

Recall that here, as in the previous sections, by an operator we mean  a bounded linear operator acting between complex Banach spaces.

%%  Def 6.1
\begin{defn} \cite{BTsk, HR, Speck1} The operators $T: X\rightarrow \widetilde{X}$ and $S: Y\rightarrow \widetilde{Y}$ are said to be \emph {(algebraically and topologically) equivalent} if and only if $T=ESF$ where $E,F$ are invertible operators. More generally, $T$ and $S$ are \emph{equivalent after extension} if and only if there exist (possibly trivial) Banach spaces $X_0$, $Y_0$, called \emph{extension spaces}, and invertible bounded linear operators $E:\widetilde{Y}\oplus Y_0\rightarrow \widetilde{X}\oplus X_0$ and  $F:X\oplus X_0\rightarrow Y\oplus Y_0$, such that
\beq \label{III.1}
\left(
  \begin{array}{cc}
    T & 0 \\
    0 & I_{X_0} \\
  \end{array}
\right)=E \left(
  \begin{array}{cc}
    S & 0 \\
    0 & I_{Y_0} \\
  \end{array}
\right) F.
\eeq

In this case we say that $T\simast S$.
\end{defn}

The relation $\simast$ is an equivalence relation. 
Operators that are equivalent after extension have many features in common. In particular, using the notation
$X\simeq Y$ to say that two Banach spaces $X$ and $Y$ are isomorphic, i.e., that there exists an invertible operator from $X$ onto $Y$, and the notation $\imag A$ to denote the range of an operator $A$, we have the following.

%%%   Theor 6.2
\begin{thm} \label{thm:3.2}
\cite{BTsk} Let $T: X\rightarrow \widetilde{X}$, $S: Y\rightarrow \widetilde{Y}$ be operators and assume that $T \simast S$. Then
\begin{enumerate}
%1
  \item $\ker T \simeq \ker S$;
 %2
  \item $\imag T$ is closed if and only if $\imag S$ is closed and, in that case,  $\widetilde{X}/\imag T \simeq \widetilde{Y}/\imag S$;
 %3
  \item if one of the operators $T$, $S$ is generalised (left, right) invertible, then the other is generalised (left, right) invertible too;
 %4 
  \item $T$ is Fredholm if and only if $S$ is Fredholm and in that case $\dim \ker T=\dim \ker S$, $\codim \imag T=\codim \imag S$.
\end{enumerate}
\end{thm}
%%%%%%%%%%%%%%%%
More properties can be found in \cite{BTsk,Speck1}, for instance.
\medskip

Now let us consider the operator
$A_g^{\alpha, \theta}: K_\theta\rightarrow K_\alpha$ and the operator
\beq \label{III.3}
P_\alpha g P_\theta+ Q_\theta : H_p^+\rightarrow K_\alpha \oplus \theta H_p^+.
\eeq
It is easy to see that
\beq \label{III.4}
A_g^{\alpha, \theta} \simast P_\alpha g P_\theta + Q_\theta
\eeq
because
\beq \label{III.5}
\left(
  \begin{array}{cc}
    A_g^{\alpha, \theta} & 0 \\
    0 & I_{\theta H_p^+} \\
  \end{array}
\right)=E_1 \left(
  \begin{array}{cc}
    P_\alpha g P_\theta+ Q_\theta & 0 \\
    0 & I_{\{0\}} \\
  \end{array}
\right) F_1
\eeq
where
\beq \label{III.6}
F_1 : K_\theta \oplus \theta H_p^+\rightarrow H_p^+ \oplus \{0\}
\eeq
\beq \label{III.7}
E_1 : (K_\alpha \oplus \theta H_p^+)\oplus \{0\}\rightarrow K_\alpha\oplus  \theta H_p^+
\eeq
are invertible operators (defined in an obvious way). On the other hand, it is clear that
\beq \label{III.8}
P_\alpha g P_\theta+ Q_\theta \simast \left(
                                          \begin{array}{cc}
                                            P_\alpha g P_\theta+ Q_\theta & 0 \\
                                            0 & P^+ \\
                                          \end{array}
                                        \right)
\eeq
where the operator on the right-hand side is defined from $(H_p^+)^2$ into $(K_\alpha \oplus \theta H_p^+) \times H_p^+$.
Now, from \eqref{III.3}
%\beq \label{III.9}
%P_\alpha g P_\theta+ Q_\theta : H_p^+\rightarrow K_\alpha \oplus \theta H_p^+.
%\eeq
we have
\beq \label{III.10}
P_\alpha g P_\theta+ Q_\theta = (P^+-P_\alpha T_g Q_\theta)(P_\alpha T_g+Q_\theta)
\eeq
where:
%% Lemma 6.3
\begin{lem} \label{lemIII.3}The operator
\beq \label{III.11}
P^+-P_\alpha T_g Q_\theta: K_\alpha \oplus \theta H_p^+ \rightarrow K_\alpha \oplus \theta H_p^+
\eeq
is invertible.
\end{lem}
%%%%%%%%%%%
\beginpf
First we prove that $P^+ \pm P_\alpha T_g Q_\theta$ maps $K_\alpha \oplus \theta H_p^+$ into $K_\alpha \oplus \theta H_p^+$.
Indeed, let $\phi_\alpha \in K_\alpha$, $\phi_+ \in H_p^+$; then
\[(P^+ \pm P_\alpha T_g Q_\theta)(\phi_\alpha+\theta \phi_+)=\phi_\alpha+\theta \phi_+ \pm P_\alpha T_g(\theta \phi_+)\]
because $Q_\theta \phi_\alpha=0$. For the same reason ($Q_\theta P_\alpha=0$), we have
\[(P^+ \pm P_\alpha T_g Q_\theta)(P^+ \mp P_\alpha T_g Q_\theta)=P^+ \mp P_\alpha T_g Q_\theta \pm P_\alpha T_g Q_\theta =P^+\]
and therefore the operator \eqref{III.11} is invertible, with inverse
\[P^+ + P_\alpha T_g Q_\theta: K_\alpha \oplus \theta H_p^+ \rightarrow K_\alpha \oplus \theta H_p^+ . \]
\endpf
%%%%%%%%%%%
Thus, with 
\[
T= \left(
 \begin{array}{cc}
  P^+-P_\alpha T_g Q_\theta & 0 \\
   0 & P^+ \\
   \end{array}
   \right),
\]
we can write
\[
\left(
 \begin{array}{cc}
  P_\alpha g P_\theta+Q_\theta & 0 \\
   0 & P^+ \\
   \end{array}
   \right)= T \left(
 \begin{array}{cc}
  P_\alpha T_g+Q_\theta & 0 \\
   0 & P^+ \\
   \end{array}
   \right)
\]

\[
=T\left(
 \begin{array}{cc}
  T_\theta & P_\alpha \\
   -P^+ & T_{\bar{\alpha}} \\
   \end{array}
   \right) \left(
 \begin{array}{cc}
  T_{\bar{\theta}} & 0 \\
   T_g-Q_\alpha(T_g-T_{\alpha \bar{\theta}}) & T_\alpha \\
   \end{array}
   \right)
\]

\beq \label{III.14}
=T\left(
 \begin{array}{cc}
  T_\theta & P_\alpha \\
   -P^+ & T_{\bar{\alpha}} \\
   \end{array}
   \right) \left(
 \begin{array}{cc}
  T_{\bar{\theta}} & 0 \\
   T_g & T_\alpha \\
   \end{array}
   \right) \left(
 \begin{array}{cc}
  P^+ & 0 \\
   -T_{\bar{\alpha}}(T_g-T_{\alpha \bar{\theta}}) & P^+ \\
   \end{array}
   \right).
\eeq

On the right-hand side of the last equality,

\;(i) the first factor, $T$, is invertible in $(K_\alpha\oplus \theta H_p^+)\times H_p^+$ by Lemma \ref{lemIII.3},

\;(ii) the second factor is invertible as an operator from $(H_p^+)^2$ into $(K_\alpha\oplus \theta H_p^+)\times H_p^+$ by Lemma  \ref{lemIII.4} below,

\;(iii) the last factor is invertible in $(H_p^+)^2$ by Lemma \ref{lemIII.5} below.

%%%  Lemma 6.4
\begin{lem} \label{lemIII.4}
The operator $T_1:(H_p^+)^2\rightarrow (K_\alpha\oplus \theta H_p^+)\times H_p^+$  defined by
\beq \label{III.15}
T_1(\phi_{1+}, \phi_{2+})=\left(
                            \begin{array}{cc}
                              T_\theta & P_\alpha \\
                              -P^+ & T_{\bar{\alpha}} \\
                            \end{array}
                          \right)\left(
                                   \begin{array}{c}
                                     \phi_{1+} \\
                                     \phi_{2+} \\
                                   \end{array}
                                 \right)
\eeq
is invertible.
\end{lem}
%%%%%
\beginpf
Given any $(\psi_{1+}, \psi_{2+}) \in (K_\alpha\oplus \theta H_p^+)\times H_p^+$, it follows from \eqref{III.15} that
\beq \label{III.16}
T_1(\phi_{1+}, \phi_{2+})=(\psi_{1+}, \psi_{2+})   
\eeq

\beq \label{III.17}
 \iff \left\{\begin{array}{c}
                          \theta \phi_{1+}+ P_\alpha \phi_{2+} = \psi_{1+} \\
                             \\
                          - \phi_{1+}+T_{\bar{\alpha}} \phi_{2+} = \psi_{2+}
                        \end{array}
 \right..
\eeq
The first equation in \eqref{III.17} implies that
\beq \label{III.18}
\theta \phi_{1+}= Q_\theta \psi_{1+}, \quad P_\alpha \phi_{2+}= P_\alpha \psi_{1+}
\eeq
and from the second equation in \eqref{III.17} we have

\beq \label{III.19}
\phi_{1+}+ \psi_{2+}= \bar{\alpha} Q_\alpha \phi_{2+};
\eeq
therefore
\beq \label{III.20}
Q_\alpha \phi_{2+}= \alpha \phi_{1+} + \alpha \psi_{2+}= \alpha \bar{\theta} Q_\theta \psi_{1+} + \alpha \psi_{2+}.
\eeq
From \eqref{III.18} and \eqref{III.20} we see that \eqref{III.16} implies that
\beq\label{III.21}
\phi_{1+}=\bar{\theta} Q_\theta \psi_{1+}, \quad  \phi_{2+}=(P_\alpha+ \alpha\bar{\theta} Q_\theta )\psi_{1+}+T_\alpha \psi_{2+}.
\eeq
It follows that $T_1$ is injective (replacing $\psi_{1+}$ and $\psi_{2+}$ by 0) and surjective (since for any $\psi_{1+}\in K_\alpha \oplus \theta H_p^+$ and any $\psi_{2+}\in H_p^+$ there exist $\phi_{1+}, \phi_{2+} \in H_p^+ $, given by \eqref{III.21}, such that \eqref{III.16} holds.

Moreover \eqref{III.21} yields an expression for the inverse operator:
\[
T_1^{-1}:(K_\alpha\oplus \theta H_p^+)\times H_p^+\rightarrow (H_p^+)^2
\]
\beq \label{III.22}
T_1^{-1}\left(
          \begin{array}{c}
            \psi_{1+} \\
            \psi_{2+} \\
          \end{array}
        \right)= \left(
                   \begin{array}{cc}
                     T_{\bar{\theta}} & 0 \\
                    P_\alpha+\alpha \bar{\theta} Q_\theta  & T_\alpha \\
                   \end{array}
                 \right)\left(
          \begin{array}{c}
            \psi_{1+} \\
            \psi_{2+} \\
          \end{array}
        \right).
\eeq
\endpf

%%%%%%    Lemma 6.5
\begin{lem}\label{lemIII.5} The operator
\beq \label{III.23}T_2:(H_p^+)^2\rightarrow (H_p^+)^2, \quad T_2=\left(
                   \begin{array}{cc}
                    P^+ & 0 \\
                    -T_{\bar{\alpha}}(T_g-T_{\alpha \bar{\theta}} ) & P^+ \\
                   \end{array}
                 \right) \eeq
is invertible, with inverse given by
 \beq \label{III.24}
 T_2^{-1}=\left(
                   \begin{array}{cc}
                    P^+ & 0 \\
                    T_{\bar{\alpha}}(T_g-T_{\alpha \bar{\theta}})  & P^+ \\
                   \end{array}
                 \right) \eeq

\end{lem}

%%%%%%
\beginpf
This follows from the fact that $T_2$ is of the form
\[\left(
    \begin{array}{cc}
      P^+ & 0 \\
      A & P^+ \\
    \end{array}
  \right)
\]
where $A$ is an operator in $H_p^+$ which commutes with $P^+$.
\endpf
%%%%%
From \eqref{III.4}, \eqref{III.8}, \eqref{III.14} and Lemmas \ref{lemIII.3}, \ref{lemIII.4} and \ref{lemIII.5} we now conclude the following.

%% Theor 6.6
\begin{thm} 
\label{thm:3.6}
$A_g^{\alpha, \theta} \simast T_G$ where $G=\left(
                                                              \begin{array}{cc}
                                                                \bar{\theta} & 0 \\
                                                                g & \alpha \\
                                                              \end{array}
                                                            \right)
.$
\end{thm}

As an immediate consequence of Theorem \ref{thm:3.6}, one may study properties of ATTO (or TTO), such as Fredholmness and invertibility, from known results for Toeplitz operators with matricial symbols and vice-versa. As an illustration, we consider the following class of TTO. Let $\theta(\xi)=e^{i\xi}, e_\lambda(\xi)=e^{i\lambda\xi}$ for $\lambda\in\mathbb R$, and
\[
g_\lambda = b\,e_{-\beta}-\lambda+ \sum_{k=1}^n (a_k \,e_{k\alpha})
\]
where $\alpha, \beta \in (0,1), \;\alpha+\beta>1,\;\alpha/\beta \notin \mathbb Q, \; b, \lambda, a_k \in \mathbb C$ for $k=1,...,n$, and $n=[1/\alpha]$ is the integer part of $1/\alpha$. By Theorem \ref{thm:3.6}, $A_{g_\lambda}^{\theta}$ is invertible, or Fredholm, if and only if the same holds for $T_{G_\lambda}$ with 
\[
G_\lambda=\left(
                                                              \begin{array}{cc}
                                                                e_{-1} & 0 \\
                                                                g_\lambda & e_1 \\
                                                              \end{array}
                                                            \right)
.
\]
For $\lambda \neq 0$, $T_{G_\lambda}$ is invertible by Theorem 5.1 in \cite{CMS}, which moreover provides explicit formulas for a bounded canonical factorisation $G_\lambda=(G_\lambda)_-\,(G_\lambda)_+$. This yields an explicit expression for the inverse operator 
\[
(T_{G_\lambda})^{-1}=(G_\lambda)_+^{-1} P^+ (G_\lambda)_-^{-1} I:(H_p^+)^2\rightarrow (H_p^+)^2, 
\]
and the inverse operator $(A_{g_\lambda}^{\theta})^{-1}$ can be obtained using \eqref{III.14}  and the related results in this section.

For $\lambda=0$ we have $G_\lambda H_+=H_-$ with $H_\pm \in (H_\infty^\pm)^2$ given by
\[
H_+=(e_\beta,\, -e_{\alpha+\beta-1} \,\sum_{k=1}^n (a_k \,e_{(k-1)\alpha}))
\]
\[
H_-=(e_{\beta -1},\,b) .
\]

By Theorem 5.3 in \cite{CDR} it follows that $\dim \ker T_{G_0}=\infty$ and therefore $T_{G_0}$ (and, consequently, $A_{g_0}^{\theta}$ ) is not Fredholm. 

Since $A_{g_\lambda}^{\theta}=A_{g_0-\lambda}^{\theta}$, we conclude that 
\[
\sigma_{ess} (A_{g_0}^{\theta})=\sigma (A_{g_0}^{\theta})=\sigma_p (A_{g_0}^{\theta})=\{0\}.
\]
%%%%%%
%%%%%%  Section 7  -  Kernels of ATTO with analytic symbols

\section{Kernels of ATTO with analytic symbols and invariant subspaces}
\label{sec:7}

\bigskip

TTO have generated much interest, and so have T-kernels 
(kernels of Toeplitz operators) -- see, for example \cite{CP14,sarason} and the 
references therein. We
are therefore led to consider kernels of ATTO. If we do so, we immediately see that, given an inner function $\theta$ and any inner function $\alpha$ such that $\alpha \preceq \theta$, we have
\beq \label{II.6}
\ker A_g^\theta \subset \ker A_g^{\alpha, \theta}
\eeq
(see Figure~1).
\begin{figure}
\begin{center}
 %%%   Fig. 1  %%%%%%%%%%%%%%%%
 \includegraphics[scale=0.3,angle=0]{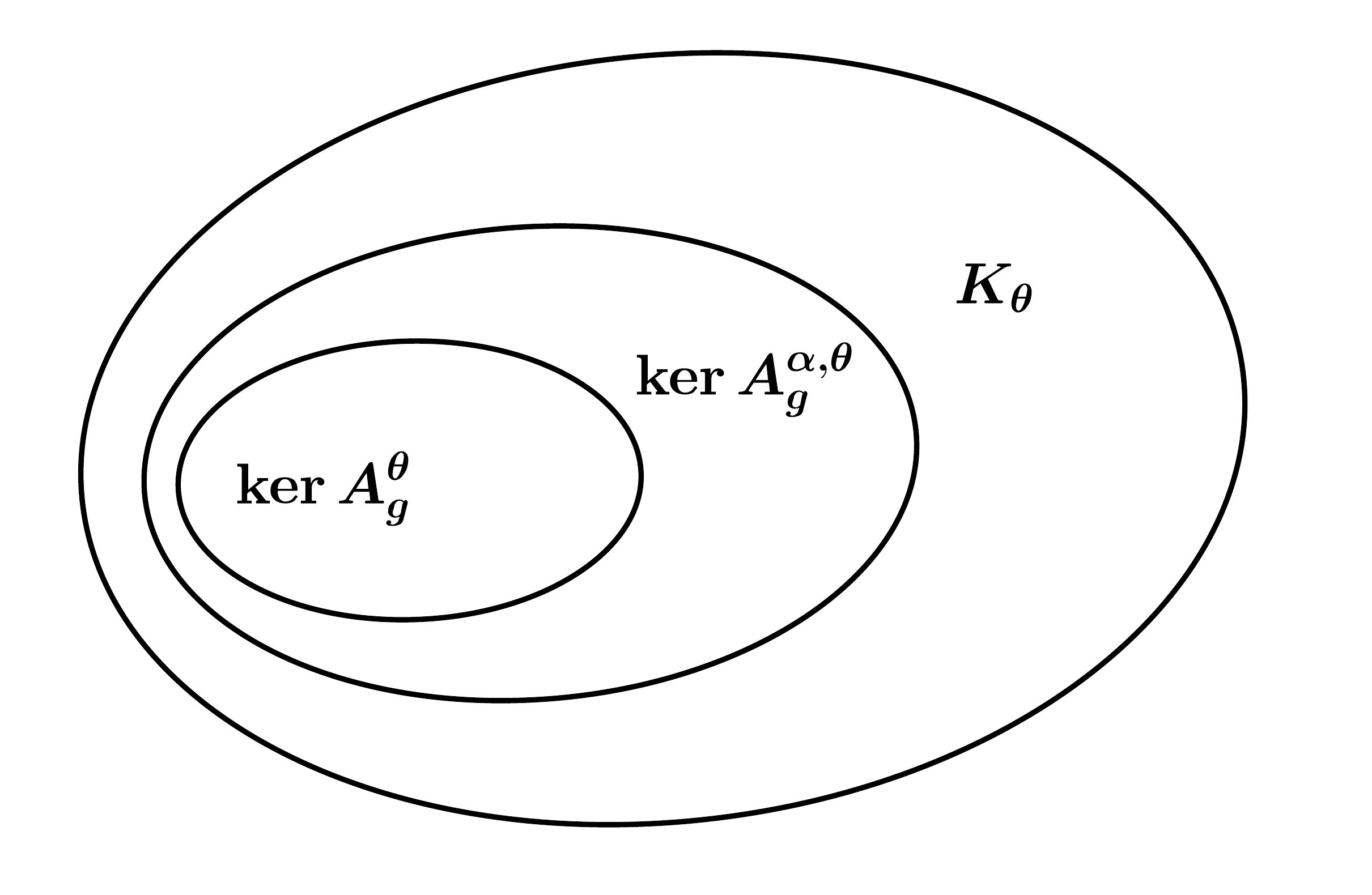}
\centerline{Fig. 1} 
\end{center}
\end{figure}

More precisely,

\beq\label{II.7}
\ker A_g^\theta =\ker A_g^{\alpha, \theta}\cap \ker B_g^{\alpha, \theta}
\eeq

where all the spaces involved are kernels of ATTO of different kinds (considering that the TTO $A_g^\theta$
is a particular case of an ATTO).

Since, according to \eqref{II.7}, $\ker A_g^{\alpha, \theta}$ is ``bigger"  than $\ker A_g^{\theta}$, it is natural to think that it may be simpler to characterize. Thus, determining the former can be seen as a first step towards determining the latter;  the elements
$\phi_+ \in \ker A_g^{\theta}$ may then be singled out by adding the condition
\[
B_g^{\alpha, \theta}\,\phi_+=0.
\]
This line of reasoning was used in \cite{CKS} to study Toeplitz operators with $2 \times 2$ triangular matrix symbols with almost periodic entries.

\bigskip

By Theorems \ref{thm:3.6} and \ref{thm:3.2},
$\ker A^{\alpha,\theta}_{g} \simeq  \ker T_G$
where $g\in L_\infty$ and
\beq\label{5.1b}
G=\left(
                                                              \begin{array}{cc}
                                                                \bar{\theta} & 0 \\
                                                                g & \alpha \\
                                                              \end{array}
                                                            \right).
\eeq
Denoting by $P_j$ the projection defined by
\[
P_j(\phi_1,\phi_2)=\phi_j \qquad (j=1,2),
\]
we have $\ker T_G  \simeq P_1(\ker T_G)$.
Indeed, $\phi_+=(\phi_{1+},\phi_{2+}) \in \ker T_G$ if and only if
we have
\[
G\phi_+ = \phi_- \qquad \hbox{with} \quad \phi_- \in (H_p^-)^2,
\]
which is equivalent to
\begin{eqnarray}
\overline\theta \phi_{1+} &=&  \phi_{1-}, \nonumber \\
g \phi_{1+} + \alpha \phi_{2+} &=& \phi_{2-},\label{eq:4.6}
\end{eqnarray}
and it is clear from (\ref{eq:4.6}) that $\phi_{1+}$   uniquely defines
$\phi_{1-}, \phi_{2+}$ and $\phi_{2-}$, since we have
\[
\phi_{1-}=\overline\theta \phi_{1+}\,, \qquad
\phi_{2-} = P^{-}(g \, \phi_{1+})=0\,, \qquad \hbox{and} \quad
\phi_{2+} = -\overline\alpha  (g \, \phi_{1+}).
\]

It is also easy to see that 
\beq\label{5.1a}
\phi_{1+} \in \ker A^{\alpha,\theta}_{g_+} \iff \phi_{1+} \in P_1(\ker T_G),
\eeq
i.e., the elements of $\ker A^{\alpha,\theta}_{g_+}$ are the
first components of the elements of $\ker T_G$, where $G$ is given by \eqref{5.1b}.

Let us now consider asymmetric truncated Toeplitz operators with symbols in $H^+_\infty$, of the form
$A^{\alpha,\theta}_{g_+}$, where
$\alpha$ and $\theta$ are inner functions such that $\alpha \preceq \theta$ and $g_+ \in H_\infty^+$.

In what follows recall that $K_{\alpha,\theta}=\alpha K_{\overline\alpha \theta}$, the shifted
model space that is the image of the projection $P_{\alpha,\theta}=P_\theta-P_\alpha$, and that
\[
K_{\theta}= K_\alpha \oplus K_{\alpha,\theta} \,,
\]
where the sum is orthogonal if $p=2$.  
The next theorem shows that shifted model spaces are the kernels of
ATTO with analytic symbols. First, however, we prove an auxiliary result.

%%  Lemma 7.1
\begin{lem}\label{L}
Given $g_+\in H_\infty^+ \setminus \{0\}$ and an inner function $\theta$,
\[
g_+ \phi_+\in \theta H_p^+ \Leftrightarrow \phi_+ \in \theta \bar \beta H_p^+ 
\]
with $\beta=GCD (g_+^i, \theta)$, where $g_+^i$ is the inner factor of the inner-outer factorization $g_+=g_+^i\, g_+^o$.
\end{lem}

\beginpf
Let $g_+\phi_+=\theta\psi_+$ with $\psi_+\in H_p^+$. Using the superscripts $i$ and $o$ to denote the inner and outer factors respectively, we have
\[
g_+^i g_+^o \phi_+^i \phi_+^o =-\theta \psi_+^i \psi_+^o,
\]
so that $g_+^i \phi_+^i = C \theta \psi_+^i$ for some $C \in \CC$ with $|C|=1$.
Dividing both sides of this equation by $\beta=\GCD(\theta,g_+^i)$ we obtain
\[
\frac{g_+^i}{\beta}\phi^i_+ = C \frac{\theta}{\beta} \psi_+^i
\]
and since $g_+^i/\beta$ and $\frac{\theta}{\beta}$ are relatively prime, it follows that $\frac{\theta}{\beta}$ 
divides $\phi_+^i$; thus $\phi_+\in \theta \bar \beta H_p^+$. Conversely, if $\phi_+= \theta \bar \beta \psi_+$ with $\psi_+\in H_p^+$, then $g_+\phi_+=(g_+^i \bar \beta)g_+^o \theta \psi_+ \in \theta H_p^+$.
\endpf

%%  Theorem 7.2
\begin{thm}
\label{thm:4.1}
Let $\alpha$ and $\theta$ be inner functions with $\alpha \preceq \theta$, and suppose that $g_+ \in H_\infty^+\setminus \{0\}$.
Then $\ker A_{g_+}^{\alpha,\theta}= K_{\gamma,\theta}$, 
with $\gamma=\alpha/\beta$ where, denoting by $g^i_+$ the inner factor in an inner-outer factorization of $g_+$,
we have $\beta=\GCD(\alpha,g_+^i)$.
\end{thm}

\beginpf
We have $\phi_{1+} \in \ker A_{g_+}^{\alpha,\theta}$ if and only if
\beq\label{eq:4.11}
\left(
                                                              \begin{array}{cc}
                                                                \bar{\theta} & 0 \\
                                                                g_+ & \alpha \\
                                                              \end{array}
                                                            \right)
\begin{pmatrix}
\phi_{1+} \\ \phi_{2+}
\end{pmatrix}
= 
\begin{pmatrix}
\phi_{1-} \\ \phi_{2-}
\end{pmatrix},
\eeq
where as usual $\phi_j^\pm \in H_p^\pm$ for $j=1,2$.
Thus $g_+ \phi_{1+} + \alpha \phi_{2+} = \phi_{2-} = 0$, and therefore
$g_+ \phi_{1+} =-\alpha \phi_{2+}$. By Lemma \ref {L}, we have $\phi_{1+} \in \gamma H_p^+$ and thus $\phi_{1+} \in \gamma H_p^+ \cap K_\theta=K_{\gamma,\theta}$\,.\\
Conversely, if $\phi_{1+} \in K_{\gamma,\theta} \subset \gamma H_p^+$\,, then by Lemma \ref {L} we have $g_+\phi_{1+}\in\alpha H_p^+$, so that we can write
\[
g_+ \phi_{1+} + \alpha \phi_{2+} = \phi_{2-}
\]
with $\phi_{2+}\in H_p^+$ and $\phi_{2-}=0$. Hence (\ref{eq:4.11})  is satisfied and 
$\phi_{1+} \in \ker A_{g_+}^{\alpha,\theta}$.
\endpf

%% Corollary 7.3
\begin{cor}
\label{cor:4.2}
Let $\alpha$ and $\theta$ be inner functions with $\alpha \preceq \theta$. Then $K_{\alpha,\theta}=\ker A_1^{\alpha,\theta}$ and $K_{\theta}=\ker A_\alpha^{\alpha,\theta}$.
\end{cor}

%% Corollary 7.4
\begin{cor}
\label{cor:4.3}
With the same assumptions as in Theorem \ref{thm:4.1}, if $p=2$ we have
\[
\ker A_{g_+}^{\alpha,\theta}=K_\theta \ominus K_\gamma = \gamma H_2^+ \ominus \theta H_2^+.
\]
\end{cor}

This holds, in particular for the TTO $A_{g_+}^\theta$, where $\alpha=\theta$, in which case we have (\cite{nik1})
\[
\ker A_{g_+}^{\theta} = \frac{\theta}{\beta} H_2^+ \ominus \theta H_2^+.
\].

Moreover, for all $p\in (1,\infty)$:

%% Corollary 7.5
\begin{cor}
\label{cor:4.4}
With the same assumptions as in Theorem \ref{thm:4.1} we have the following:

(i) $A_{g_+}^{\alpha,\theta}=0$ if and only if $g_+\in \alpha H_\infty^+$;

(ii) $A_{g_+}^{\alpha,\theta}$ is injective if and only if $\alpha=\theta$ and $\beta$ is a constant;

(iii) $\dim \ker A_{g_+}^{\alpha,\theta}<\infty$ if and only if $\overline \alpha \theta$ and $\beta$ are finite Blaschke products and, in that case, $\dim \ker A_{g_+}^{\alpha,\theta}=n_1+n_2$ where $n_1$ and $n_2$ are the number of zeroes of $\overline \alpha \theta$ and $\beta$, respectively.

(iv) For $\alpha=\theta$, $\dim \ker A_{g_+}^{\theta}<\infty$ if and only if $\beta$ is a finite Blaschke product and, in that case, $\dim \ker A_{g_+}^{\theta}$ is equal to the number of common zeroes of $g_+^i$ and $\theta$.
\end{cor}

As an immediate consequence we see that, in the particular case of the truncated shift
with symbol $r$ given by (\ref{def:r}), we have
$\ker A_r^\theta=\{0\}$ if $\theta(i)\neq 0$, and $\ker A_r^\theta=\frac{\theta}{r}K_r=\spam \{\frac {\theta}{\xi-i}\}$ if $\theta(i)=0$.

\bigskip

Shifted model spaces are also associated with ATTO in a different way: they are
the (closed) invariant subspaces of the truncated shift $A^\theta_r$.

%% Theorem 7.6
\begin{thm}
\label{thm:4.5}
The lattice $\Lat(A_r^\theta)$ consists of the
spaces $K_{\alpha,\theta}$, where $\alpha \preceq \theta$.
\end{thm}

\beginpf
For $\alpha \preceq \theta$ and $\beta=\theta \overline\alpha$, we have $K_{\alpha,\theta}=\alpha K_\beta$;  let $k^+$ be any function in $K_\beta$. Then $k^+=P_\beta \,\phi_+$ for some $\phi_+\in H_p^+$ and
\[
P_\theta\, r(\alpha k^+)=P_\theta \,r(\alpha P_\beta \,\phi_+)\\
=P_\theta\, r P_\theta \,\alpha \phi_+=P_\theta \,r\alpha \phi_+\\
=\alpha P_\beta (r\phi_+)\in \alpha K_\beta.
\]
Thus every space $K_{\alpha,\theta}$ is invariant for $A_r^\theta$.
To show the converse, we begin with the observation that for the Hardy space $H_p(\DD)$ of the unit disc, we have a version of
Beurling's theorem for each $1<p<\infty$; namely that the nontrivial invariant subspaces of the shift $T_z$ are all of the form
$\alpha H_p$ for some inner function $\alpha$. See, for example, \cite[Cor. C.2.1.20]{nik}. By means of the standard
isometric isomorphism between $H_p(\DD)$ and $H_p^+$ 
given in (\ref{eq:Viso})
we see that the
same result holds for the shift $T_r$ on $H_p^+$.

Next, using the duality between $H_p^+$ and $H_q^+$ (up to isomorphism), we see that the $T_r^*$-invariant
subspaces in $H_q^+$ are the annihilators of the invariant subspaces for $T_r$, i.e., the
model spaces
\[
K^q_\alpha = \{f \in H_q^+: \int_\RR f \overline g = 0 \ \forall g \in \alpha H^p\} = \alpha H_q^- \cap H_q^+.
\]

Now if $A_r^\theta$ is a restricted shift on $H_p^+$, then its Banach space adjoint is the restriction of $T_r^*$
to its invariant subspace $K_\theta^q$, so that its adjoint has invariant subspaces $K_\alpha^q$ where $\alpha \preceq \theta$.

Using duality once more we conclude that the invariant subspaces of $A_r^\theta$ take the form
\[
\{ f \in K_\theta^p: \int_\RR f\overline g=0 \ \forall g \in K_\alpha^q \} = K^p_\theta \cap \alpha H^p = K_{\alpha,\theta}\,,
\]
where   $\alpha \preceq \theta$.

\endpf

%% Corollary 7.7
\begin{cor}
$\Lat(A_r^\theta)=\{\ker A_{g_+}^{\alpha,\theta} :\alpha\preceq\theta,\, g_+\in H_\infty^+\}$.
\end{cor}

We may now prove a theorem of Lax--Beurling flavour for the ``truncated shift'' semigroup 
on $K_\theta$ given by
\[
T(t)=A^\theta_{e_t}, \qquad (t \ge 0),
\]
where $e_t \in H_\infty^+$ is the 
inner function given by $e_t(\xi)=e^{it\xi}$.

\begin{thm}
The common invariant subspaces of the semigroup $(T(t))_{t \ge 0}$ are the 
shifted model spaces $K_{\alpha,\theta}$, where $\alpha \preceq \theta$.
\end{thm}

\beginpf
It is easy to see that these subspaces are all invariant under the semigroup,
since if $\alpha$ divides a function
$f \in K_\theta$
then it also divides $ T(t)f$.\\

The converse is proved as in \cite[Thm.~3.1.5]{lols}, the standard Lax--Beurling theorem.
By writing
% is=\xi
\[
\frac{1}{\xi+i} = \frac{1}{i}\int_0^\infty e^{-t}e^{i t\xi} \, dt,
\]
and approximating the integral by Riemann sums, we see that the ATTO operator with symbol $1/(\xi+i)$
is the strong limit of a sequence of finite linear combinations of the ATTO with symbols $e_t$.
Hence any closed subspace invariant under the semigroup is also invariant under $A^\theta_r$, and thus
is a shifted model space, as required.
\endpf

\section*{Acknowledgments}

This work was partially supported by Funda\c{c}\~{a}o para a Ci\^{e}ncia e a Tecnologia (FCT/Portugal), through Project
PTDC/MAT/121837/2010 and PEst-OE/EEI/LA0009/2013.

%This work was partially supported by Funda\c{c}\~{a}o para a Ci\^{e}ncia e a Tecnologia (FCT/Portugal), through Project
%PTDC/MAT/121837/2010 and  Project Est-C/MAT/UI0013/2011, and by the Research Centre of Mathematics of the University of Minho through the FEDER Funds �Programa Operacional Factores de Competitividade COMPETE�.

\end{document}